\documentclass{amsart}
\usepackage{graphicx}
\newtheorem{thm}{Theorem}[section]
\newtheorem{cor}[thm]{Corollary}
\newtheorem{lem}[thm]{Lemma}
\newtheorem{prop}[thm]{Proposition}
\theoremstyle{definition}

\theoremstyle{remark}
\newtheorem{rem}[thm]{Remark}
\numberwithin{equation}{section}
\begin{document}

\title[Optimal control of a large dam]{Optimal control of a large dam,
taking into account the water costs}%
\author{Vyacheslav M. Abramov}%
\address{School of Mathematical Sciences, Monash University, Building 28M, Wellington road, Clayton, VIC 3800, Australia}%
\email{vyacheslav.abramov@sci.monash.edu.au}%

%\thanks{}%

%\date{}%
%\dedicatory{}%
%\commby{}%
% ----------------------------------------------------------------
\begin{abstract}
Consider a dam model, $L^{upper}$ and $L^{lower}$ are upper and,
respectively, lower levels, $L = L^{upper}-L^{lower}$ is large and
if  the level of water is between these bounds, then the dam is
said to be in a normal state.
 Passage across lower or upper levels leads to damage. Let $J_1=j_1L$ and $J_2=j_2L$
denote the damage costs per time unit of crossing the lower and,
correspondingly, upper level where $j_1$ and $j_2$ are given real
constants. It is assumed that input stream of water is described
by a Poisson process, while the output stream is state dependent.
Let $L_t$ denote the level of water in time $t$, and $c_{L_t}$
denote the water cost at level $L_t$ ($L^{lower}<L_t\leq
L^{upper}$). Assuming that
$p_1=\lim_{t\to\infty}\mathbf{P}\{L_t=L^{lower}\}$,
$p_2=\lim_{t\to\infty}\mathbf{P}\{L_t>L^{upper}\}$ and
$q_i=\lim_{t\to\infty}\mathbf{P}\{L_t=i\}$ ($L^{lower}<i\leq
L^{upper}$) exist, the aim of the paper is to choose the
parameters of an output stream (specifically defined in the paper)
minimizing the long-run expenses
$$J=p_1J_1+p_2J_2+\sum_{i=L^{lower}+1}^{L^{upper}}q_ic_i.$$
% The more particular problem, not taking
%into account the water costs, has been recently studied in
%[Abramov, \emph{J. Appl. Prob.}, 44 (2007), 249-258]. The present
%paper partially answers the question: \textit{How does the
%structure of water costs affect the optimal solution?}

\end{abstract}
\maketitle
% ----------------------------------------------------------------
\section{Introduction}\label{Introduction}
A large dam is defined by the parameters $L^{lower}$ and
$L^{upper}$, which are, respectively, the lower and upper levels
of the dam. If the current level is between these bounds, the dam
is assumed to be in a normal state. The difference $L =
L^{upper}-L^{lower}$ is large, and this is the reason for calling
the dam \textit{large}. This property enables us to use asymptotic
analysis as $L\to\infty$ and solve easier different problems of
optimal control than we would were than the dam is not large.
%one of such simplest problems of optimal control for a
%large dam has been solved in \cite{Abramov 2007}.

Let $L_t$ denote the water level in time $t$. If
$L^{lower}<L_t\leq L^{upper}$, then the state of the dam is called
\textit{normal}. Passage across lower or upper level leads to
damage. The costs per time unit of this damage is $J_1=j_1L$ and
$J_2=j_2L$ for lower and upper levels correspondingly, where $j_1$
and $j_2$ are given real constants.
 The water inflow is described by the Poisson process with rate
 $\lambda$. In practice, this means that the arrival of water units
 is registered by counter at random instants $t_1$, $t_2$, \ldots, and the
 times between consecutive instants are mutually independent
and exponentially distributed with parameter $\lambda$.

 The outflow of water is state-dependent as
 follows. If the level of water is between $L^{lower}$ and
 $L^{upper}$, then an interval between departures of units of
 water (inverse output flow) has the probability distribution
 function
 $B_1(x)$. If level of water exceeds $L^{upper}$, then an
 inverse output flow has the probability distribution function $B_2(x)$. The probability
 distribution function $B_2(x)$ is assumed to obey the condition
 $\int_0^\infty x\mbox{d}B_2(x)<\frac{1}{\lambda}$. If
the level of water is $L^{lower}$ exactly, then output of water is
 frozen and it resumes again as soon as the level of water exceeds
 the level $L^{lower}$.

Let $c_{L_t}$ denote the cost of water at level $L_t$. The
sequence $c_i$ is assumed to be positive and non-increasing. The
problem of the present paper is to choose the parameter
$\int_0^\infty x\mbox{d}B_1(x)$ of the dam in the normal state
minimizing the objective function
\begin{equation}\label{I1}
J = p_1J_1 + p_2J_2 + \sum_{i=L^{lower}+1}^{L^{upper}}c_iq_i,
\end{equation}
where
\begin{eqnarray}
p_1&=&\lim_{t\to\infty}\mathbf{P}\{L_t=L^{lower}\},\label{I2}\\
p_2&=&\lim_{t\to\infty}\mathbf{P}\{L_t>L^{upper}\},\label{I3}\\
q_i&=&\lim_{t\to\infty}\mathbf{P}\{L_t=L^{lower}+i\}, \
i=1,2,\ldots,L\label{I4}.
\end{eqnarray}
In the queueing formulation, the level $L^{lower}$ is identified
with an empty queue, and the dam model is the following queueing
system with service depending on queue-length. If immediately
before a service begins the queue-length exceeds the level $L$,
then the customer is served by the probability distribution
$B_2(x)$. Otherwise, the service time distribution is $B_1(x)$.
The value $p_1$ is the stationary probability of empty system, the
value $p_2$ is the stationary probability that a customer is
served by probability distribution $B_2(x)$, and $q_i$,
$i=1,2,\ldots,L$, are the stationary probabilities of the
queue-length process, so $p_1+p_2+\sum_{i=1}^L q_i=1$. (For the
described queueing system, the right-hand side limits in relations
\eqref{I2}-\eqref{I4} do exist.)

In our study, the parameter $L$ increases indefinitely, and we
deal with the series of queueing systems. The parameters above,
such as $p_1$, $p_2$, $J_1$, $J_1$ as well as other parameters are
functions of $L$. The argument $L$ will be often omitted in these
functions.

Similarly to \cite{Abramov 2007}, it is assumed that the input
parameter $\lambda$ and probability distribution function $B_2(x)$
are given, while the appropriate probability function $B_1(x)$
should be chosen from the specified parametric family of functions
$B_1(x,C)$. (Actually, we deal with the family of probability
distributions $B_1(x)$ depending on two parameters $\delta$ and
$L$ in series. Then the parametric family of distributions
$B_1(x,C)$ is described by the family of possible limits of
$\delta L$ as $\delta\to0$ and $L\to\infty$.)

 The outflow rate
associated with the probability distribution function $B_1(x)$
should be chosen such that the minimum of the objective function
of \eqref{I1} is associated with the choice of the parameter $C$
resulting in the choice of the corresponding probability
distribution function $B_1(x,C)$.

The more particular problem, where the objective function has the
form $J= p_1J_1 + p_2J_2$ (i.e. the water costs are not taken into
account), has been studied in Abramov \cite{Abramov 2007}. In this
case the solution to the control problem is unique and has one of
the following three forms. Denote $\rho_2=\lambda\int_0^\infty
x\mbox{d}B_2(x)$ and $\rho_1=\rho_1(C)=\lambda\int_0^\infty
x\mbox{d}B_1(x,C)$. Then, in the case
$j_1=j_2\frac{\rho_2}{1-\rho_2}$, the optimal solution is
$\rho_1=1$. In the case $j_1>j_2\frac{\rho_2}{1-\rho_2}$, the
optimal solution has the form $\rho_1=1+\delta$, where $\delta(L)$
is a small positive parameter, and $\delta(L)L\to C$ as
$L\to\infty$. In the case $j_1<j_2\frac{\rho_2}{1-\rho_2}$, the
optimal strategy has the form $\rho_1=1-\delta$, and
$\delta(L)L\to C$ as $L\to\infty$. The parameter $C$ is a unique
solution of a specific minimization problem precisely formulated
in \cite{Abramov 2007}. It has been also shown in \cite{Abramov
2007} that the solution to the control problem is insensitive to
the type of probability distributions $B_1(x)$ and $B_2(x)$.
Specifically, it is expressed via the first moment of $B_2(x)$ and
the first two moments of $B_1(x)$.

Following \cite{Abramov 2007}, we use the notation
$\rho_{1,l}=\lambda^l\int_0^\infty x^l\mbox{d}B_1(x)$, $l=$2,3.
The existence of a moment of the order corresponding to
$\rho_{1,l}$ will be specially assumed in formulations of
statements corresponding to case studies.

The problem studied in the present paper substantially
distinguishes from that studied in \cite{Abramov 2007}. Although
the both control problem of the present paper and \cite{Abramov
2007} are closely related, the new components of the present
problem change the problem substantially. The problem of the
present paper requires a much deepen and delicate analysis (for
example, it is demonstrated in the next section that the
asymptotic methods of the earlier paper \cite{Abramov 2007} do not
longer work here, and one should use more delicate techniques
instead), and the main results on optimal control policies are
deepen as well.

Essential difficulty of the control problem in the present
formulation is to prove a uniqueness of the optimal solution,
while in the case of the particular problem of \cite{Abramov
2007}, the uniqueness of the solution follows automatically from
the explicit representations of the functionals obtained there.

It is assumed in the present paper that $c_i$ is a non-increasing
sequence. If the cost sequence $c_i$ were an arbitrary bounded
sequence, then a richer class of possible cases could be studied.
However, in the case of arbitrary cost sequence, the solution need
not be unique, and arbitrary costs $c_i$, say increasing in $i$,
seem not to be useful and, therefore, are not considered here. The
practical applications of the results obtained in the paper are
not restricted by the area of water research. The results of the
present paper are meaningful for many specific inventory and
storage problems as well. Moreover, the model considered in this
paper describes more realistically a production/realization
process of large warehouses rather than an inflow/outflow process
of large dams, because in many cases of a water inflow the
seasonality is important. Nevertheless, similarly to \cite{Abramov
2007}, the problem formulation in this paper is given in terms of
water research.

More realistic models arising in practice assume that the
probability distribution function $B_1(x)$ should also depend on
$i$, i.e have representation $B_{1,i}(x)$. The model of the
present paper, where $B_1(x)$ is the same for all $i$, under
appropriate additional information can approximate those more
general models. Namely, one can suppose the stationary service
time distribution $B_1(x)$ has the representation
$B_1(x)=\sum_{i=1}^Lq_iB_{1,i}(x)$ ($q_i$, $i=1,2,\ldots,L$ are
the state probabilities), and the solution to the control problem
for $B_1(x)$ enables us to find then the approximate solutions to
the control problem for $B_{1,i}(x)$, $i=1,2,\ldots,L$ by using
the Bayes rule.

Similarly to the solution of the control problem of \cite{Abramov
2007}, the solution of the present problem with extended criteria
\eqref{I1} is related to the same class of solutions as in
\cite{Abramov 2007}. That is, it must be either $\rho_1=1$ or one
of two limits of $\rho_1=1+\delta$, $\rho_1=1-\delta$ for positive
small vanishing $\delta$ as $L$ increases indefinitely, and
$L\delta\to C$. The reason for this is that the penalties upon
reaching upper or lower level are of order $O(L)$ (i.e. increase
to infinity as $L\to\infty$), while the water costs are assumed to
be bounded as $L$ tends to infinity, and although the water costs
affect the solution of the control problem, this influence remains
in the framework of the same class of solutions mentioned above.

The following new questions are of special interest here.

\smallskip
1. What is the structure of an optimal solution? Is an optimal
solution unique?

\smallskip
The answer to these questions is the main result of the paper. The
questions are answered by Theorem \ref{thm3}. We prove that a
solution to the control problem does exist and unique, however
there are some additional mild assumptions related to the class of
probability distributions $\{B_1(x)\}$. The proof of the existence
and uniqueness of a solution is based on special techniques of
Mathematical Analysis. Specifically, we use the known techniques
of majorization inequalities \cite{Hardy Littlewood Polya 1952},
\cite{Marschall Olkin 1979} in order to prove the monotonicity of
specified functions. This property of monotonicity is then used to
prove a uniqueness of a solution.

\smallskip

2. Under what relation between $j_1$, $j_2$, $\rho_2$, $c_i$ (and
maybe other parameters of the model) the optimal strategy is
$\rho_1=1$?

\smallskip
If the water costs are not taken into account, then the condition
for $\rho_1=1$ is $j_1=j_2\frac{\rho_2}{1-\rho_2}$. This result
has been proved in \cite{Abramov 2007}. It has a simple intuitive
explanation and is a consequence of the well-known property of the
stream of lost calls during a busy period of M/GI/1/$n$ queues,
under the assumption that the expected interarrival and service
times are equal (see Abramov \cite{Abramov 1997} as well as
Righter \cite{Righter 1999} or Wolff \cite{Wolff 2002}). Taking
into account the structure of water costs generally changes this
condition for the aforementioned optimal solution $\rho_1=1$. We
prove that the optimal solution $\rho_1=1$ is achieved under the
condition $j_1\leq j_2\frac{\rho_2}{1-\rho_2}$, and the equality
in this relation holds if and only if the water costs are the same
at all levels of water. This result is only the partial answer to
the question. More exact answers can be obtained in particular
cases, and one of them is the case of linearly decreasing costs as
the level of water increases (for brevity, this case is called
\textit{linear costs}). In the case of linear costs we derive more
exact and useful representations, which enable us to calculate
numerically the relation between $j_1$ and $j_2$ to have finally
the optimal solution $\rho_1=1$. The relevant numerical results
are provided for special values of the parameters of the model.

The rest of the paper is organized as follows. In Section
\ref{Stationary probabilities}, the asymptotic behavior of the
stationary probabilities is studied. This section is structured
into four subsections. In Section \ref{Preliminaries}, some
necessary results related to asymptotic properties of
characteristic of the state-dependent queueing system are recalled
for modelling the behavior of a large dam. In Section \ref{Case
1}, elementary asymptotic properties of the state probabilities in
the case $\rho_1=1$ are established. The analysis of this case is
based on application of a Tauberian theorem of Postnikov
\cite{Postnikov 1980} for the convolution type recurrence
relation. In Section \ref{Case 2} the behavior of the state
probabilities is studied for the case $\rho_1=1+\delta$, where
positive parameter $\delta$ vanishes as $L$ tends to infinity. The
analysis of this case uses known results on asymptotic behavior of
characteristics of the state-dependent queueing system obtained in
Section \ref{Preliminaries}. Delicate asymptotic analysis is
provided in Section \ref{Case 3} for the case $\rho_1=1-\delta$,
where positive parameter $\delta$ vanishes as $L$ tends to
infinity. For the analysis of this case we involve the asymptotic
result of Willmot \cite{Willmot 1988}, which is also used in the
present paper. All of the aforementioned asymptotic results of
Section \ref{Stationary probabilities} are used to establish
asymptotic properties of special functionals. Then we solve the
problem of minimization of these functionals. In Section
\ref{Solution}, we solve the control problem. We prove existence
and uniqueness of the solution under mild assumptions on the class
of probability distributions $\{B_1(x)\}$. We establish a
structure of the solution and completely answer the question 1
posed in the introduction. In section \ref{Examples}, the case of
linear costs is studied. We establish explicit representations for
functionals of the control problem and provide numerical answer to
the question 2 posed in the introduction.

\section{Stationary probabilities of the state-dependent queueing
system and their asymptotic behavior}\label{Stationary
probabilities}

In this section, the explicit expressions are derived for the
stationary probabilities, and their asymptotic behavior is
studied. These results will be used in our further findings of the
optimal solution.

\subsection{Preliminaries}\label{Preliminaries} We first recall some
results of \cite{Abramov 2007} and then develop them in order to
obtain the explicit representations for the stationary
probabilities that will be required in our further analysis. Let
$T_L$, $\nu_L$, $I_L$ denote correspondingly a busy period, the
number of served customers during a busy period, and an idle
period. Let $\nu_L^{(1)}$ and $\nu_L^{(2)}$ denote the number of
customers served during a busy period by probability distribution
functions $B_1(x)$ and $B_2(x)$ correspondingly, and let
$T_L^{(1)}$ and $T_L^{(2)}$ denote the time spent for the service
of customers during a busy period by the probability distribution
functions $B_1(x)$ and $B_2(x)$ correspondingly.

It was shown in \cite{Abramov 2007} that the probabilities $p_1$
and $p_2$ can be expressed explicitly via $\mathbf{E}\nu_L^{(1)}$,
$\rho_1$ and $\rho_2$ only. Namely,
\begin{equation}\label{SP1}
p_1=\frac{1-\rho_2}{1+(\rho_1-\rho_2)\mathbf{E}\nu_L^{(1)}},
\end{equation}
and
\begin{equation}\label{SP2}
p_2=\frac{\rho_2+\rho_2(\rho_1-1)\mathbf{E}\nu_L^{(1)}}{1+(\rho_1-\rho_2)\mathbf{E}\nu_L^{(1)}}.
\end{equation}

Such kind of representations is convenient, because
$\mathbf{E}\nu_L^{(1)}$ satisfies the convolution type recurrence
relation
\begin{equation}\label{SP3}
\mathbf{E}\nu_L^{(1)}=\sum_{j=0}^L\mathbf{E}\nu_{L-j+1}^{(1)}\int_0^\infty\mbox{e}^{-\lambda
x}\frac{(\lambda x)^j}{j!}\mbox{d}B_1(x), \ \
\mathbf{E}\nu_0^{(1)}=1,
\end{equation}
where $\mathbf{E}\nu_n^{(1)}$ denotes the expectation of the
number of served customers during a busy period of the M/GI/1/$n$
queue ($n$=0,1,\ldots). Recurrence relation \eqref{SP3} is in turn
a special case of the recurrence relations
\begin{equation}\label{SP4}
Q_n=\sum_{j=0}^nQ_{n-j+1}r_j,
\end{equation}
with $r_0>0$, $r_j\geq0$ for all $j\geq1$, and $r_0+r_1+\ldots=1$.
The detailed theory of these recurrence relations can be found in
Tak\'{a}cs \cite{Takacs 1967}. Asymptotic behavior of $Q_n$ as
$n\to\infty$ has been studied by Tak\'{a}cs \cite{Takacs 1967} and
Postnikov \cite{Postnikov 1980}.

The stationary probabilities $q_i$ can be obtained from the
renewal arguments (e.g. Ross \cite{Ross 2000}). Namely, for
$i=1,2,\ldots,L$ we have
\begin{equation}\label{SP5}
q_i=\frac{\mathbf{E}T_i^{(1)}-\mathbf{E}T_{i-1}^{(1)}}{\mathbf{E}T_L+\mathbf{E}I_L},
\end{equation}
where $\mathbf{E}T_i^{(1)}$ denotes the expectation of a busy
period of the M/GI/1/$i$ queue ($i$=0,1,\ldots). The probabilities
given by \eqref{SP5} can be also written
\begin{equation}\label{SP6}
q_i=\rho_1\frac{\mathbf{E}\nu_i^{(1)}-\mathbf{E}\nu_{i-1}^{(1)}}{\mathbf{E}\nu_L}.
\end{equation}
Indeed, from Wald's equation \cite{Feller 1966}, p.384 we have
$\mathbf{E}T_i^{(k)}$ =
$\frac{\rho_k}{\lambda}\mathbf{E}\nu_i^{(k)}$, $k=1,2$, and
therefore the numerator of \eqref{SP5} is rewritten
\begin{equation}\label{SP7}
\mathbf{E}T_i^{(1)}-\mathbf{E}T_{i-1}^{(1)}=\frac{\rho_1}{\lambda}
\left(\mathbf{E}\nu_i^{(1)}-\mathbf{E}\nu_{i-1}^{(1)}\right).
\end{equation}
Next, using the fact that the number of arrivals during a busy
cycle coincides with the number of served customers during a busy
period, from Wald's equation \cite{Feller 1966}, p.384 we obtain
\begin{equation}\label{SP8}
\lambda\mathbf{E}T_L+\lambda\mathbf{E}I_L=\mathbf{E}\nu_L.
\end{equation}
Therefore, \eqref{SP6} is the consequence of \eqref{SP7} and
\eqref{SP8}.

On the other hand, \eqref{SP8} can be rewritten
\begin{equation}\label{SP9}
\begin{aligned}
\lambda\mathbf{E}T_L+\lambda\mathbf{E}I_L&=\lambda\mathbf{E}T_L+1\\
&=\lambda\left(\mathbf{E}T_L^{(1)}+\mathbf{E}T_L^{(2)}\right)\\
&=\rho_1\mathbf{E}\nu_L^{(1)}+\rho_{2}\mathbf{E}\nu_L^{(2)}+1.
\end{aligned}
\end{equation}
From \eqref{SP8} and \eqref{SP9} we have the equation
\begin{equation*}\label{SP10}
\mathbf{E}\nu_L^{(1)}+\mathbf{E}\nu_L^{(2)}=\rho_1\mathbf{E}\nu_L^{(1)}+\rho_{2}\mathbf{E}\nu_L^{(2)}+1,
\end{equation*}
resulting in
\begin{equation}\label{SP11}
\mathbf{E}\nu_L^{(2)}=\frac{1}{1-\rho_2}-\frac{1-\rho_1}{1-\rho_2}\mathbf{E}\nu_L^{(1)}
\end{equation}
and, consequently by adding to the both sides of \eqref{SP11}
$\mathbf{E}\nu_L^{(1)}$, we obtain
\begin{equation}\label{SP12}
\mathbf{E}\nu_L=\frac{1}{1-\rho_2}+\frac{\rho_1-\rho_2}{1-\rho_2}\mathbf{E}\nu_L^{(1)}.
\end{equation}
Substituting \eqref{SP12} for \eqref{SP6}, we finally obtain
\begin{equation}\label{SP13}
q_i=\frac{\rho_1(1-\rho_2)}{1+(\rho_1-\rho_2)\mathbf{E}\nu_L^{(1)}}
\left(\mathbf{E}\nu_i^{(1)}-\mathbf{E}\nu_{i-1}^{(1)}\right), \ \
i=1,2,\ldots,L.
\end{equation}
Comparison with \eqref{SP1} enables us to rewrite \eqref{SP13} in
the other form
\begin{equation}\label{SP14}
q_i=\rho_1p_1
\left(\mathbf{E}\nu_i^{(1)}-\mathbf{E}\nu_{i-1}^{(1)}\right), \ \
i=1,2,\ldots,L
\end{equation}
as well.

In order to study the asymptotic behavior of the loss
probabilities, let us recall the earlier results on asymptotic
behavior of $\mathbf{E}\nu_L^{(1)}$ as $L$ increases indefinitely
(e.g. \cite{Abramov 1997}, \cite{Abramov 2004}, \cite{Abramov
2007}). This asymptotic behavior can be obtained from the results
of Tak\'acs \cite{Takacs 1967}, p.22-23, on asymptotic behavior of
recurrence relation \eqref{SP4} as $n\to\infty$. (The formulation
of this Tak\'acs theorem can be found in \cite{Abramov 2007} as
well.)

If $\rho_1=1$ and $\rho_{1,2}<\infty$, then
\begin{equation}\label{SP16}
\lim_{L\to\infty}\frac{\mathbf{E}\nu_L^{(1)}}{L}=\frac{2}{\rho_{1,2}},
\end{equation}
and if $\rho_1>1$, then
\begin{equation}\label{SP17}
\lim_{L\to\infty}\left[\mathbf{E}\nu_L^{(1)}-\frac{1}{\varphi^L(1+\lambda\widehat
B_1^\prime(\lambda-\lambda\varphi))}\right]=\frac{1}{1-\rho_1},
\end{equation}
where $\widehat B_1(s)$ is the Laplace-Stieljes transform of
$B_1(x)$, $\varphi$ is the least in absolute value root of the
functional equation $z=\widehat B_1(\lambda-\lambda z)$, and
$\widehat B_1^\prime(s)$ denotes the derivative of $\widehat
B_1(s)$.

We do not consider the case $\rho_1<1$, because it is not
meaningful for our analysis. Only relations \eqref{SP16} and
\eqref{SP17} are useful for the asymptotic analysis of the state
probabilities, which is provided below.

\subsection{The case $\mathbf{\rho_1=1}$}\label{Case 1}
Asymptotic relation \eqref{SP16} is not enough in order to obtain
an asymptotic behavior of the state probabilities. As it was
mentioned above, this asymptotic relation \eqref{SP16} is obtained
from the aforementioned Tak\'acs theorem \cite{Takacs 1967}, which
has been also used in \cite{Abramov 2007} for the analysis of the
asymptotic behavior of $p_1$ and $p_2$. The asymptotic analysis of
the state probabilities is more delicate than that analysis of
$p_1$ and $p_2$, and asymptotic relation \eqref{SP16} is not
enough for this purpose. Therefore, we have to use the more
precise asymptotic estimation given by the Tauberian theorem of
Postnikov \cite{Postnikov 1980} on the asymptotic behavior of
recurrence sequence \eqref{SP4}.

\begin{lem}\label{lem2}(Postnikov \cite{Postnikov 1980}, Sec. 25.)
Let $\gamma_1=\sum_{n=1}^\infty nr_n=1$,
$\gamma_2=\sum_{n=2}^\infty n(n-1)r_n<\infty$, and $r_0+r_1<1$.
Let $Q_0\neq0$ be an arbitrarily given initial value of recurrence
sequence \eqref{SP4}. Then as $n\to\infty$
\begin{equation*}\label{SP18}
Q_{n+1}-Q_n=\frac{2Q_0}{\gamma_2}+o(1).
\end{equation*}
%(It is used the notation: $\gamma_i=\sum_{j=i}^\infty\prod_{l=1}^i
%(l+j-i)r_j$.)
\end{lem}

From this lemma we obtain the following statement.

\begin{cor}\label{cor1}
 For any $j\geq0$,
\begin{equation}\label{SP19}
\mathbf{E}\nu_{L-j}^{(1)}-\mathbf{E}\nu_{L-j-1}^{(1)}=\frac{2}{\rho_{1,2}}+o(1),
\ \text{as} \ L\to\infty.
\end{equation}
\end{cor}

\begin{proof}
In the case where $r_j=\int_0^\infty\mbox{e}^{-\lambda
x}\frac{(\lambda x)^j}{j!}\mbox{d}B_1(x)$, $j=0,1,\ldots$, the
Tauberian condition $r_0+r_1<1$ of Lemma \ref{lem2} is always
fulfilled, and the result of corollary follows directly from Lemma
\ref{lem2}. The proof of this fact uses the theory of analytic
functions and has been provided in \cite{Abramov 2004} (Theorem
4.6) and \cite{Abramov 2002} (Theorem 3.3). Since this proof is
elementary and short it is repeated here again.

We must prove that for \textit{some} $\lambda_0>0$ the equality
\begin{equation}\label{SP19-add}
\int_0^\infty\mbox{e}^{-\lambda_0 x}\left(1+\lambda_0
x\right)\mbox{d}B_1(x)=1
\end{equation}
is not the case. If it is so, then only the inequality
\begin{equation*}
\int_0^\infty\mbox{e}^{-\lambda x}\left(1+\lambda
x\right)\mbox{d}B_1(x)<1
\end{equation*}
must hold.

Indeed, $\int_0^\infty\mbox{e}^{-\lambda x}\left(1+\lambda
x\right)\mbox{d}B_1(x)$ is an analytic function in $\lambda$, and
therefore, according to the theorem on the maximum module of an
analytic function, equality \eqref{SP19-add} must hold for
\textit{all} $\lambda_0\geq0$. This means, that \eqref{SP19-add}
is valid if and only if
$\int_0^\infty\mbox{e}^{-\lambda_0x}\frac{(\lambda_0x)^j}{j!}\mbox{d}B_1(x)$=0
for all $j\geq2$ and $\lambda_0\geq0$. In this case the
Laplace-Stieltjes transform $\widehat{B}_1(\lambda)$ must be a
linear function in $\lambda$, i.e.
$\widehat{B}_1(\lambda)=d_0+d_1\lambda$, $d_0$ and $d_1$ are some
constants. However, since $|\widehat{B}_1(\lambda)|\leq1$, we have
$d_0=1$ and $d_1=0$. This is a trivial case where $B_1(x)$ is
concentrated in point 0, and therefore it is not a probability
distribution function having a positive mean. Thus
\eqref{SP19-add} is not the case, and the statement of the
corollary follows.
\end{proof}
%
%Relation \eqref{SP19} can be also rewritten in the more convenient
%form as
%\begin{equation}\label{SP19+}
%\mathbf{E}\nu_{j+1}^{(1)}-\mathbf{E}\nu_{j}^{(1)}=\frac{2}{\rho_{1,2}}+o(1),
%\end{equation}
%{as} $L\to\infty$ ($j=1,2,\ldots,L-1$), which is used in the
%sequel.

From \eqref{SP19} and \eqref{SP1}, \eqref{SP14} we have the
following lemma.
\begin{lem}\label{lem3}In the case $\rho_1=1$ and $\rho_2<\infty$
for any $j\geq0$ we have
\begin{equation}\label{SP20}
\lim_{L\to\infty}Lq_{L-j}=1.
\end{equation}
\end{lem}

\begin{proof} Asymptotic relation \eqref{SP20} follows immediately from \eqref{SP1},
\eqref{SP14} and \eqref{SP19}.
\end{proof}

Combining the result obtained in \cite{Abramov 2007} and Lemma
\ref{lem3} we have the asymptotic result for the functional $J$ =
$J_1p_1$ + $J_2p_2$ + $\sum_{i=1}^L q_ic_i$.
\begin{prop}\label{prop1}
In the case $\rho_1=1$ and $\rho_2<\infty$ we have
\begin{equation}\label{SP21}
\lim_{L\to\infty}J(L)=j_1\frac{\rho_{1,2}}{2}+j_2\frac{\rho_2}{1-\rho_2}\cdot\frac{\rho_{1,2}}{2}+c^*,
\end{equation}
where
\begin{equation*}\label{SP21-add}
c^*=\lim_{L\to\infty}\frac{1}{L}\sum_{i=1}^L c_i.
\end{equation*}
\end{prop}

\begin{proof}The first two term of the right-hand side of \eqref{SP21} easily
follow from \eqref{SP1}, \eqref{SP2} and asymptotic representation
\eqref{SP16}. The last term of the right-hand side of \eqref{SP21}
follows from Lemma \ref{lem3}, since
\begin{equation*}
\lim_{L\to\infty}\sum_{i=1}^L
q_ic_i=\lim_{L\to\infty}\frac{1}{L}\sum_{i=1}^L c_i=c^*.
\end{equation*}
Notice, that the representation for the first two terms of the
right-hand side of \eqref{SP21} has a relation to the particular
case of the dam model not taking into account the water costs and
been obtained in \cite{Abramov 2007}.
\end{proof}

\begin{rem}
The statement of Proposition \ref{prop1} confirms that the class
of the possible solutions of the control problem is the same as
that for the particular problem studied in \cite{Abramov 2007},
since according to Proposition \ref{prop1} the optimal value of
the functional $J$ is finite (i.e. does not increase to infinity
as $L$ increases indefinitely) and cannot exceed the right-hand
side of \eqref{SP21}.
\end{rem}

\subsection{The case $\mathbf{\rho_1=1+\delta,
\delta>0}$}\label{Case 2} In the case $\rho_1=1+\delta$ we have
the following.

\begin{thm}\label{thm1}
Assume that $\rho_1=1+\delta$, $\delta>0$, and $L\delta\to C>0$ as
$\delta\to 0$ and $L\to\infty$. Assume that $\rho_{1,3}(L)$ is a
bounded sequence, and there exists
$\widetilde{\rho}_{1,2}=\lim_{L\to\infty}\rho_{1,2}(L)$. Then, for
any $j\ge0$
\begin{equation}\label{SP22}
q_{L-j}=\frac{\mbox{e}^{2C/\widetilde{\rho}_{1,2}}}{\mbox{e}^{2C/\widetilde{\rho}_{1,2}}-1}
\left(1-\frac{2\delta}{\widetilde{\rho}_{1,2}}\right)^j\frac{2\delta}{\widetilde{\rho}_{1,2}}
%\cdot\frac{\mbox{e}^{2C/\widetilde{\rho}_{1,2}}}{\mbox{e}^{2C/\widetilde{\rho}_{1,2}}-1}
+o(\delta).
\end{equation}
\end{thm}

\begin{proof}
Expanding first \eqref{SP17} for large $L$, we have
\begin{equation}\label{SP23}
\begin{aligned}
\mathbf{E}\nu_{L-j}^{(1)}=\frac{\varphi^j}{\varphi^{L}\left[1+\lambda\widehat{B}_1^\prime(\lambda-\lambda\varphi)\right]}+
\frac{1}{1-\rho_1}+o(1).
\end{aligned}
\end{equation}
From \eqref{SP23} for large $L$ we have
\begin{equation}\label{SP23-add}
\mathbf{E}\nu_{L-j}^{(1)}-\mathbf{E}\nu_{L-j-1}^{(1)}=
\frac{(1-\varphi)\varphi^j}{\varphi^L\left[1+\lambda\widehat{B}_1^\prime(\lambda-\lambda\varphi)\right]}+o(1).
\end{equation}

Next, under the condition of the theorem we have the following
expansion
\begin{equation}\label{SP24}
\varphi=1-\frac{2\delta}{\widetilde{\rho}_{1,2}}+O(\delta^2).
\end{equation}
(The details of this expansion can be found in \cite{Subhankulov
1976}, p.326 or \cite{Abramov 2006}, p.21. See also the proof of a
similar fact in Lemma \ref{lem4} below.) Then, taking into account
\eqref{SP24}, we also have
\begin{equation}\label{SP25}
1+\lambda\widehat{B}_1^\prime(\lambda-\lambda\varphi)=\delta+O(\delta^2).
\end{equation}
Next, taking into account asymptotic expansions \eqref{SP24} and
\eqref{SP25} from \eqref{SP23} and \eqref{SP23-add} we have
\begin{equation}\label{SP25-add}
\mathbf{E}\nu_L^{(1)}=\frac{\mbox{e}^{2C/\widetilde{\rho}_{1,2}}-1}{\delta}
+O(1),
\end{equation}
and for any $j=0,1,\ldots$
\begin{equation}\label{SP25-add2}
\mathbf{E}\nu_{L-j}^{(1)}-\mathbf{E}\nu_{L-j-1}^{(1)}=\mbox{e}^{2C/\widetilde{\rho}_{1,2}}
\left(1-\frac{2\delta}{\widetilde{\rho}_{1,2}}\right)^j\frac{2}{\widetilde{\rho}_{1,2}}~[1+o(1)].
\end{equation}
Now, using explicit representation \eqref{SP13}, for any
$j=0,1\ldots$ we obtain the desired asymptotic relation
\begin{equation*}
q_{L-j}=\frac{\mbox{e}^{2C/\widetilde{\rho}_{1,2}}}{\mbox{e}^{2C/\widetilde{\rho}_{1,2}}-1}
\left(1-\frac{2\delta}{\widetilde{\rho}_{1,2}}\right)^j\frac{2\delta}{\widetilde{\rho}_{1,2}}~[1+o(1)].
\end{equation*}
The theorem is proved.
\end{proof}

Using the result of \cite{Abramov 2007} we have the following.

\begin{prop}\label{prop2} Under the assumptions of Theorem \ref{thm1}
denote the objective function $J$ by $J^{upper}$. We have the
following representation
\begin{equation}\label{SP29}
\begin{aligned}
J^{upper}=&C\left[j_1\frac{1}{\mbox{e}^{2C/\widetilde{\rho}_{1,2}}-1}+j_2\frac{\rho_2\mbox{e}^{2C/\widetilde{\rho}_{1,2}}}
{(1-\rho_2)({\mbox{e}^{2C/\widetilde{\rho}_{1,2}}-1})}\right]+c^{upper},
\end{aligned}
\end{equation}
where
\begin{equation}\label{SP29+}
c^{upper}=\frac{2C}{\widetilde{\rho}_{1,2}}\cdot
\frac{\mbox{e}^{2C/\widetilde{\rho}_{1,2}}}{\mbox{e}^{2C/\widetilde{\rho}_{1,2}}-1}
\lim_{L\to\infty}\frac{1}{L}~\widehat{C}_L\left(1-\frac{2C}{\widetilde{\rho}_{1,2}L}\right),
\end{equation}
and $\widehat C_L(z)=\sum_{j=0}^{L-1} c_{L-j}z^j$ is a backward
generating cost function.
\end{prop}

\begin{proof}
The representation for the term
$$
C\left[j_1\frac{1}{\mbox{e}^{2C/\widetilde{\rho}_{1,2}}-1}+j_2\frac{\rho_2\mbox{e}^{2C/\widetilde{\rho}_{1,2}}}
{(1-\rho_2)({\mbox{e}^{2C/\widetilde{\rho}_{1,2}}-1})}\right]
$$
is known from \cite{Abramov 2007}. The new term in \eqref{SP29},
which takes into account the water costs, is $c^{upper}$.
Therefore, keeping in mind representation \eqref{SP22}, for this
term we obtain:
$$
c^{upper}=\lim_{L\to\infty}\sum_{j=0}^{L-1}
q_{L-j}c_{L-j}=\lim_{L\to\infty}\sum_{j=0}^{L-1} c_{L-j}\cdot
\frac{\mbox{e}^{2C/\widetilde{\rho}_{1,2}}}{\mbox{e}^{2C/\widetilde{\rho}_{1,2}}-1}
\left(1-\frac{2\delta
L}{\widetilde{\rho}_{1,2}L}\right)^j\frac{2\delta
L}{\widetilde{\rho}_{1,2}L},
$$
and representation \eqref{SP29+} follows.
\end{proof}

\subsection{The case $\mathbf{\rho_1=1-\delta, \delta>0}$}\label{Case 3} The study of
this case is more delicate and based on a special analysis. Our
additional assumption here is that the class of probability
distribution functions $\{B_1(x)\}$ is given such that there
exists a unique root $\tau>1$ of the equation
\begin{equation}\label{SP30}
z=\widehat{B}(\lambda-\lambda z),
\end{equation}
and there exists the first derivative
$\widehat{B}^\prime(\lambda-\lambda\tau)$.

 Under the assumption that $\rho_1<1$ the root of \eqref{SP30} is
not necessarily exists. Such type of condition has been considered
by Willmot \cite{Willmot 1988} to obtain the asymptotic behavior
for high queue-level probabilities in stationary M/GI/1 queues.
Denote the stationary probabilities of M/GI/1 queues by
$q_i[M/GI/1]$, $i=0,1,\ldots$. There was shown in \cite{Willmot
1988} that
\begin{equation}\label{SP32}
q_i[M/GI/1]=\frac{(1-\rho_1)(1-\tau)}{\tau^i[1+\lambda\widehat{B}^\prime(\lambda-\lambda\tau)]}[1+o(1)],
\ \ \text{as} \ i\to\infty.
\end{equation}
%(Both the numerator and denominator of the right-hand side of
%\eqref{SP32} are negative.)
On the other hand, according to the Pollaczek-Khintchine formula
(e.g. Tak\'acs \cite{Takacs 1962}, p. 242) $q_i[M/GI/1]$ can be
represented explicitly
\begin{equation}\label{SP33}
q_i[M/GI/1]=(1-\rho_1)\left(\mathbf{E}\nu_{i}^{(1)}-\mathbf{E}\nu_{i-1}^{(1)}\right),
\  \ i=1,2,\ldots.
\end{equation}
(Representation \eqref{SP33} can be easily checked, since it
easily follows from \eqref{SP3} that
\begin{equation}\label{SP33-add}
\sum_{j=0}^\infty\mathbf{E}\nu_{j}^{(1)}z^j=\frac{\widehat{B}_1(\lambda-\lambda
z)}{\widehat{B}_1(\lambda-\lambda z)-z},
\end{equation}
and multiplication of the right-hand side of \eqref{SP33-add} by
$(1-\rho_1)(1-z)$ leads to the well-known Pollaczek-Khintchine
formula.)

From \eqref{SP32} and \eqref{SP33} we also have the following
asymptotic proportion. For large $L$ and any $j\geq0$
\begin{equation}\label{SP34}
\frac{\mathbf{E}\nu_{L-j}^{(1)}-\mathbf{E}\nu_{L-j-1}^{(1)}}
{\mathbf{E}\nu_{L}^{(1)}-\mathbf{E}\nu_{L-1}^{(1)}}=\tau^j[1+o(1)].
\end{equation}

Now we formulate and prove a theorem on asymptotic behavior of
stationary probabilities $q_i$ in the case $\rho_1=1-\delta$. In
this theorem we assume that the class of probability distributions
$\{B_1(x)\}$ is defined according to the above convention. In the
case where $\rho_1=1-\delta$, $\delta>0$ and vanishing $\delta$ as
$L$ increases indefinitely, this means that there exists
$\epsilon_0>0$ (small enough) such that for all
$0\leq\epsilon\leq\epsilon_0$, the above family of probability
distribution functions $B_{1,\epsilon}(x)$ (depending now on
parameter $\epsilon$) satisfies the following properties. Let
$\widehat{B}_{1,\epsilon}(s)$ denote the Laplace-Stieltjes
transform of $B_{1,\epsilon}(x)$. We assume that any
$\widehat{B}_{1,\epsilon}(s)$ is an analytic function in a small
neighborhood of zero, and
\begin{equation}\label{SP34-add}
\widehat{B}_{1,\epsilon}^{\prime}(-\lambda\epsilon)<\infty.
\end{equation}
%and as $\epsilon\downarrow0$
%$$
%\widehat{B}_{1,\epsilon}(-\lambda\epsilon)=1+\epsilon+
%\lambda^2\widehat{B}_{1,\epsilon}^{\prime\prime}(0)\epsilon^2+O(\epsilon^3).
%\leqno(ii)
%$$
Property \eqref{SP34-add} is required for the existence of the
probabilities $q_i$. Relation \eqref{SP32} contains the term
$\widehat{B}^\prime(\lambda-\lambda\tau)$, and for
$\tau=1+\epsilon$ this term should be finite.
%Property (ii) is the
%result of the Taylor expansion for small $\epsilon$.
Choice of the
small parameter $\epsilon$ is closely connected with that choice
of parameter $\delta$ in the theorem below.

%
%
%
%
%where $\widehat{B}_{1,\epsilon}^{\prime}(\cdot)$ denotes the
%derivative of the Laplace-Stieltjes transform
%$\widehat{B}_{1,\epsilon}(\cdot)$. The change of parameter
%$\epsilon$ corresponds to the change of vanishing parameter
%$\delta$, so as $\delta$ vanishes, $\epsilon$ vanishes as well
%(see formulation of Lemma \ref{lem4} below) and $\delta_0$ of that
%lemma corresponds to the value $\epsilon_0$.

\begin{thm}\label{thm2}
Assume that the class of probability distribution functions
$\{B_1(x)\}$ is defined according to the above convention and
satisfies \eqref{SP34-add}, $\rho_1=1-\delta$, $\delta>0$, and
$L\delta\to C>0$, as $\delta\to 0$ and $L\to\infty$. Assume that
$\rho_{1,3}=\rho_{1,3}(L)$ is a bounded sequence, and there exists
$\widetilde{\rho}_{1,2}=\lim_{L\to\infty}\rho_{1,2}(L)$. Then,
\begin{equation}\label{SP35}
q_{L-j}=\frac{2\delta}{\widetilde{\rho}_{1,2}}\cdot\frac{1}{\mbox{e}^{2C/\widetilde{\rho}_{1,2}}-1}
\left(1+\frac{2\delta}{\widetilde{\rho}_{1,2}}\right)^j~[1+o(1)]
\end{equation}
for any $j\geq0$.
\end{thm}

\begin{proof} We start the proof from the following auxiliary
result, similar to that used to prove Theorem \ref{thm1}.

\begin{lem}\label{lem4}
Under the assumptions of Theorem \ref{thm2}, the following
asymptotic representation holds
\begin{equation}\label{SP36}
\tau=1+\frac{2\delta}{\widetilde{\rho}_{1,2}}+O(\delta^2).
\end{equation}
\end{lem}

\textit{Proof of Lemma \ref{lem4}}. The proof of this lemma is
completely similar to the proof of the result in \cite{Subhankulov
1976}, p.326 mentioned in the proof of Theorem \ref{thm1} (see
relation \eqref{SP24}. According to the above convention, the
equation $z=\widehat{B}(\lambda-\lambda z)$ has a unique solution
$\tau>1$. Clearly, the root approaches 1 as $\delta$ vanishes,
because otherwise we would have a contradiction with Corollary
\ref{cor1}. Therefore, by the Taylor expansion of this equation
around the point $z=1$, we have
\begin{equation*}\label{SP37}
z=1-(1-\delta)(1-z)+\frac{\widetilde{\rho}_{1,2}}{2}(1-z)^2+O(1-z)^3.
\end{equation*}
Ignoring the last term $O(1-z)^3$ of this expansion, we have the
quadratic equation, and the two solutions of the equation are
$z=1$ and $z=1+\frac{2\delta}{\widetilde{\rho}_{1,2}}$. Therefore,
\begin{equation*}
\tau=1+\frac{2\delta}{\widetilde{\rho}_{1,2}}+O(\delta^2),
\end{equation*}
and the lemma is therefore proved. \ \ $\Box$

Let us continue the proof of the theorem. From \eqref{SP34} and
\eqref{SP36} of Lemma \ref{lem4} we have:
\begin{equation}\label{SP38}
q_{L-j}=q_L\left(1+\frac{2\delta}{\widetilde{\rho}_{1,2}}\right)^j[1+o(1)].
\end{equation}
Taking into account that
\begin{equation*}\label{SP39}
\begin{aligned}
\sum_{j=0}^{L-1}\left(1+\frac{2\delta}{\widetilde{\rho}_{1,2}}\right)^j&=
\frac{\widetilde{\rho}_{1,2}}{2\delta}\left[\left(1+\frac{2\delta}{\widetilde{\rho}_{1,2}}\right)^L-1\right]\\
&=\frac{\widetilde{\rho}_{1,2}}{2\delta}
\left(\mbox{e}^{2C/\widetilde{\rho}_{1,2}}-1\right),
\end{aligned}
\end{equation*}
from the normalization condition $p_1+p_2+\sum_{j=1}^L q_j$=1 and
the fact that the both $p_1$ and $p_2$ have the order $O(\delta)$,
we obtain:
\begin{equation}\label{SP40}
q_L=\frac{2\delta}{\widetilde{\rho}_{1,2}}\cdot\frac{1}{\mbox{e}^{2C/\widetilde{\rho}_{1,2}}-1}~[1+o(1)].
\end{equation}
The desired statement of the theorem follows from \eqref{SP40}.
\end{proof}

Using the result of \cite{Abramov 2007} we have the following.

\begin{prop}\label{prop3}
Under the assumptions of Theorem \ref{thm2} denote the objective
function $J$ by $J^{lower}$. We have the following representation
\begin{equation}\label{SP41}
J^{lower}=C\left[j_1\mbox{e}^{\widetilde{\rho}_{1,2}/2C}+j_2\frac{\rho_2}{1-\rho_2}
\left(\mbox{e}^{\widetilde{\rho}_{1,2}/2C}-1\right)\right]+c^{lower},
\end{equation}
where
\begin{equation}\label{SP42}
c^{lower}=\frac{2C}{\widetilde{\rho}_{1,2}}\cdot\frac{1}{\mbox{e}^{2C/\widetilde{\rho}_{1,2}}-1}
\lim_{L\to\infty}\frac{1}{L}~\widehat{C}_{L}\left(1+\frac{2C}{\widetilde{\rho}_{1,2}L}\right),
\end{equation}
and $\widehat{C}(z)=\sum_{j=0}^{L-1}c_{L-j}z^j$ is a backward
generating cost function.
\end{prop}

\begin{proof}
The representation for the term
\begin{equation*}
C\left[j_1\mbox{e}^{\widetilde{\rho}_{1,2}/2C}+j_2\frac{\rho_2}{1-\rho_2}
\left(\mbox{e}^{\widetilde{\rho}_{1,2}/2C}-1\right)\right]
\end{equation*}
is known from \cite{Abramov 2007}. The new term in \eqref{SP41},
which takes into account the water costs, is $c^{lower}$. Keeping
in mind representation \eqref{SP35}, for this term we obtain:
$$
c^{lower}=\lim_{L\to\infty}\sum_{j=0}^{L-1}
q_{L-j}c_{L-j}=\lim_{L\to\infty}\sum_{j=0}^{L-1} c_{L-j}\cdot
\frac{2\delta
L}{\widetilde{\rho}_{1,2}L}\cdot\frac{1}{\mbox{e}^{2C/\widetilde{\rho}_{1,2}}-1}
\left(1+\frac{2\delta L}{\widetilde{\rho}_{1,2}L}\right)^j,
$$
and representation \eqref{SP42} follows.
\end{proof}

\section{Solution of the control problem and its
properties}\label{Solution} In this section we discuss the
solution to the control problem and study its properties. The
functionals $J^{upper}$ and $J^{lower}$ are given by \eqref{SP29}
and \eqref{SP41}. The last terms of these functionals are
\begin{equation}\label{SCP1}
c^{upper}=\frac{2C}{\widetilde{\rho}_{1,2}}\cdot
\frac{\mbox{e}^{2C/\widetilde{\rho}_{1,2}}}{\mbox{e}^{2C/\widetilde{\rho}_{1,2}}-1}
\lim_{L\to\infty}\frac{1}{L}\sum_{j=0}^{L-1}
c_{L-j}\left(1-\frac{2\delta}{\widetilde{\rho}_{1,2}}\right)^j,
\end{equation}
and
\begin{equation}\label{SCP2}
c^{lower}=\frac{2C}{\widetilde{\rho}_{1,2}}\cdot\frac{1}{\mbox{e}^{2C/\widetilde{\rho}_{1,2}}-1}
\lim_{L\to\infty}\frac{1}{L}\sum_{j=0}^{L-1}
c_{L-j}\left(1+\frac{2\delta}{\widetilde{\rho}_{1,2}}\right)^j
\end{equation}
correspondingly, where $C$ = $\lim_{L\to\infty} L\delta$.

Now we give other representations for these terms. Denote
\begin{equation}\label{SCP3}
\psi(C)=\lim_{L\to\infty}\dfrac{\sum_{j=0}^{L-1}c_{L-j}\left(1-\dfrac{2C}{\widetilde{\rho}_{1,2}L}\right)^j}
{\sum_{j=0}^{L-1}\left(1-\dfrac{2C}{\widetilde{\rho}_{1,2}L}\right)^j},
\end{equation}
and
\begin{equation}\label{SCP4}
\eta(C)=\lim_{L\to\infty}\dfrac{\sum_{j=0}^{L-1}c_{L-j}\left(1+\dfrac{2C}{\widetilde{\rho}_{1,2}L}\right)^j}
{\sum_{j=0}^{L-1}\left(1+\dfrac{2C}{\widetilde{\rho}_{1,2}L}\right)^j}.
\end{equation}
Since $\{c_i\}$ is a bounded sequence, then the both limits of
\eqref{SCP3} and \eqref{SCP4} do exist.

We have the following lemma.

\begin{lem}\label{lem5} We have
\begin{equation}\label{SCP5}
c^{upper}=\psi(C),
\end{equation}
and
\begin{equation}\label{SCP6}
c^{lower}=\eta(C).
\end{equation}
\end{lem}

\begin{proof}
From \eqref{SCP3} and \eqref{SCP4} we correspondingly have the
representations
\begin{equation}\label{SCP7}
\lim_{L\to\infty}\frac{1}{L}\sum_{j=0}^{L-1}
c_{L-j}\left(1-\frac{2C}{\widetilde{\rho}_{1,2}L}\right)^j=\psi(C)
\lim_{L\to\infty}\frac{1}{L}\sum_{j=0}^{L-1}
\left(1-\frac{2C}{\widetilde{\rho}_{1,2}L}\right)^j,
\end{equation}
and
\begin{equation}\label{SCP8}
\lim_{L\to\infty}\frac{1}{L}\sum_{j=0}^{L-1}
c_{L-j}\left(1+\frac{2C}{\widetilde{\rho}_{1,2}L}\right)^j=\eta(C)
\lim_{L\to\infty}\frac{1}{L}\sum_{j=0}^{L-1}
\left(1+\frac{2C}{\widetilde{\rho}_{1,2}L}\right)^j.
\end{equation}
The desired results follow by direct transformations of the
corresponding right-hand sides of \eqref{SCP7} and \eqref{SCP8}.

Prove \eqref{SCP5}. For the right-hand side of \eqref{SCP7} we
have
\begin{equation}\label{SCP9}
\psi(C) \lim_{L\to\infty}\frac{1}{L}\sum_{j=0}^{L-1}
\left(1-\frac{2C}{\widetilde{\rho}_{1,2}L}\right)^j=\psi(C)
\left(1-\mbox{e}^{-2C/\widetilde{\rho}_{1,2}}\right)\frac{\widetilde{\rho}_{1,2}}{2C}.
\end{equation}
On the other hand, from \eqref{SCP1} we have
\begin{equation}\label{SCP10}
c^{upper}\left(1-\mbox{e}^{-2C/\widetilde{\rho}_{1,2}}\right)
\frac{\widetilde{\rho}_{1,2}}{2C}=
\lim_{L\to\infty}\frac{1}{L}\sum_{j=0}^{L-1}
c_{L-j}\left(1-\frac{2C}{\widetilde{\rho}_{1,2}L}\right)^j.
\end{equation}
Hence, from \eqref{SCP7}, \eqref{SCP9} and \eqref{SCP10} we obtain
\eqref{SCP3}. The proof of \eqref{SCP6} is completely analogous
and uses representations \eqref{SCP2} and \eqref{SCP8}.
\end{proof}

The next lemma establishes the main properties of functions
$\psi(C)$ and $\eta(C)$.

\begin{lem}\label{lem6}
The function $\psi(C)$ is a non-increasing function, and its
maximum is $\psi(0)=c^*$. The function $\eta(C)$ is an increasing
function, and its minimum is $\eta(0)=c^*$. (Recall that
$c^*=\lim_{L\to\infty}\frac{1}{L}\sum_{i=1}^L c_i$ is defined in
Proposition \ref{prop1}.)
\end{lem}

\begin{proof} Let us first prove that $\psi(0)=c^*$ is a maximum
of $\psi(C)$. For this purpose we use the following well-known
inequality (e.g. Hardy, Littlewood and Polya \cite{Hardy
Littlewood Polya 1952} or Marschall and Olkin \cite{Marschall
Olkin 1979}). Let $\{a_n\}$ and $\{b_n\}$ be arbitrary sequences,
one of them is increasing and another decreasing. Then for any
finite sum we have
\begin{equation}\label{SCP11}
\sum_{n=1}^l a_nb_n\leq\frac{1}{l}\sum_{n=1}^l a_n\sum_{n=1}^l
b_n.
\end{equation}

Applying inequality \eqref{SCP11} to finite sums of the left-hand
side of \eqref{SCP7} and passing to limit as $L\to\infty$, we have
\begin{equation}\label{SCP12}
\begin{aligned}
&\lim_{L\to\infty}\frac{1}{L}\sum_{j=0}^{L-1}
c_{L-j}\left(1-\frac{2C}{\widetilde{\rho}_{1,2}L}\right)^j\\
&\leq \lim_{L\to\infty}\frac{1}{L}\sum_{j=0}^{L-1}
c_{L-j}\lim_{L\to\infty}\frac{1}{L}\sum_{j=0}^{L-1}
\left(1-\frac{2C}{\widetilde{\rho}_{1,2}L}\right)^j\\
&=\psi(0)\lim_{L\to\infty}\frac{1}{L}\sum_{j=0}^{L-1}
\left(1-\frac{2C}{\widetilde{\rho}_{1,2}L}\right)^j.
\end{aligned}
\end{equation}
Then, comparing \eqref{SCP7} with \eqref{SCP12} enables us to
conclude,
\begin{equation*}\label{SCP13}
\psi(0)=c^*\geq\psi(C),
\end{equation*}
i.e. $\psi(0)=c^*$ is the maximum value of $\psi(C)$.

Prove now, that $\psi(C)$ is a non-increasing function, i.e. for
any nonnegative $C_1\leq C$ we have $\psi(C)\leq\psi(C_1)$.

To prove this note, that for small positive $\delta_1$ and
$\delta_2$ we have (1-$\delta_1$-$\delta_2$) = (1-$\delta_1$)
(1-$\delta_2$) + $O(\delta_1\delta_2)$. Using this idea, one can
prove the monotonicity of $\psi(C)$ by replacing
\begin{equation*}
\left(1-\frac{2C}{\widetilde{\rho}_{1,2}L}\right)=
\left(1-\frac{2C_1}{\widetilde{\rho}_{1,2}L}\right)
\left(1-\frac{2C-2C_1}{\widetilde{\rho}_{1,2}L}\right)+O\left(\frac{1}{L^2}\right),
\ \ C>C_1
\end{equation*}
in the above asymptotic relations for large $L$. Indeed, notice
that
\begin{equation}\label{SCP14.1}
\begin{aligned}
&\lim_{L\to\infty}\frac{1}{L}\sum_{j=0}^{L-1}
\left(1-\frac{2C_1}{\widetilde{\rho}_{1,2}L}\right)^j
\left(1-\frac{2C-2C_1}{\widetilde{\rho}_{1,2}L}\right)^j\\
&= \lim_{L\to\infty}\frac{1}{L}\sum_{j=0}^{L-1}
\left(1-\frac{2C_1}{\widetilde{\rho}_{1,2}L}\right)^j
\lim_{L\to\infty}\frac{1}{L}\sum_{j=0}^{L-1}
\left(1-\frac{2C-2C_1}{\widetilde{\rho}_{1,2}L}\right)^j
\end{aligned}
\end{equation}
Therefore, for any non-decreasing sequence $a_j$
\begin{equation}\label{SCP14.2}
\begin{aligned}
&\lim_{L\to\infty}\frac{1}{L}\sum_{j=0}^{L-1}
a_j\left(1-\frac{2C_1}{\widetilde{\rho}_{1,2}L}\right)^j
\left(1-\frac{2C-2C_1}{\widetilde{\rho}_{1,2}L}\right)^j\\
&\leq \lim_{L\to\infty}\frac{1}{L}\sum_{j=0}^{L-1}a_j
\left(1-\frac{2C_1}{\widetilde{\rho}_{1,2}L}\right)^j
\lim_{L\to\infty}\frac{1}{L}\sum_{j=0}^{L-1}
\left(1-\frac{2C-2C_1}{\widetilde{\rho}_{1,2}L}\right)^j
\end{aligned}
\end{equation}
Indeed, assume for contrary that
\begin{equation}\label{SCP14.2+}
\begin{aligned}
&\lim_{L\to\infty}\frac{1}{L}\sum_{j=0}^{L-1}
a_j\left(1-\frac{2C_1}{\widetilde{\rho}_{1,2}L}\right)^j
\left(1-\frac{2C-2C_1}{\widetilde{\rho}_{1,2}L}\right)^j\\
&> \lim_{L\to\infty}\frac{1}{L}\sum_{j=0}^{L-1}a_j
\left(1-\frac{2C_1}{\widetilde{\rho}_{1,2}L}\right)^j
\lim_{L\to\infty}\frac{1}{L}\sum_{j=0}^{L-1}
\left(1-\frac{2C-2C_1}{\widetilde{\rho}_{1,2}L}\right)^j
\end{aligned}
\end{equation}
Then, applying the inequality \eqref{SCP11} to the right-hand side
of \eqref{SCP14.2+}, we obtain:
\begin{equation}\label{SCP14.2++}
\begin{aligned}
&\lim_{L\to\infty}\frac{1}{L}\sum_{j=0}^{L-1}a_j
\left(1-\frac{2C_1}{\widetilde{\rho}_{1,2}L}\right)^j
\lim_{L\to\infty}\frac{1}{L}\sum_{j=0}^{L-1}
\left(1-\frac{2C-2C_1}{\widetilde{\rho}_{1,2}L}\right)^j\\
&\leq\lim_{L\to\infty}\frac{1}{L}\sum_{j=0}^{L-1}a_j\lim_{L\to\infty}\frac{1}{L}\sum_{j=0}^{L-1}
\left(1-\frac{2C_1}{\widetilde{\rho}_{1,2}L}\right)^j
\lim_{L\to\infty}\frac{1}{L}\sum_{j=0}^{L-1}
\left(1-\frac{2C-2C_1}{\widetilde{\rho}_{1,2}L}\right)^j\\
&=\lim_{L\to\infty}\frac{1}{L}\sum_{j=0}^{L-1}a_j\lim_{L\to\infty}\frac{1}{L}\sum_{j=0}^{L-1}
\left(1-\frac{2C}{\widetilde{\rho}_{1,2}L}\right)^j.
\end{aligned}
\end{equation}
Comparison of the last obtained term with the left-hand side of
\eqref{SCP14.2+} enables us to write:
\begin{equation*}\label{SCP14.2+++}
\begin{aligned}
&\lim_{L\to\infty}\frac{1}{L}\sum_{j=0}^{L-1}
a_j\left(1-\frac{2C}{\widetilde{\rho}_{1,2}L}\right)^j
>\lim_{L\to\infty}\frac{1}{L}\sum_{j=0}^{L-1}a_j\lim_{L\to\infty}\frac{1}{L}\sum_{j=0}^{L-1}
\left(1-\frac{2C}{\widetilde{\rho}_{1,2}L}\right)^j.
\end{aligned}
\end{equation*}
The contradiction with the basic inequality \eqref{SCP11} proves
\eqref{SCP14.2}.

Taking into account \eqref{SCP14.1} and \eqref{SCP14.2} the
extended version of \eqref{SCP12} after application \eqref{SCP11}
now looks
\begin{equation}\label{SCP14}
\begin{aligned}
&\lim_{L\to\infty}\frac{1}{L}\sum_{j=0}^{L-1}
c_{L-j}\left(1-\frac{2C}{\widetilde{\rho}_{1,2}L}\right)^j\\
&=\lim_{L\to\infty}\frac{1}{L}\sum_{j=0}^{L-1}
c_{L-j}\left(1-\frac{2C_1}{\widetilde{\rho}_{1,2}L}\right)^j
\left(1-\frac{2C-2C_1}{\widetilde{\rho}_{1,2}L}\right)^j\\
&\leq \lim_{L\to\infty}\frac{1}{L}\sum_{j=0}^{L-1}
c_{L-j}\left(1-\frac{2C_1}{\widetilde{\rho}_{1,2}L}\right)^j
\lim_{L\to\infty}\frac{1}{L}\sum_{j=0}^{L-1}
\left(1-\frac{2C-2C_1}{\widetilde{\rho}_{1,2}L}\right)^j\\
&=\psi(C_1)\lim_{L\to\infty}\frac{1}{L}\sum_{j=0}^{L-1}
\left(1-\frac{2C_1}{\widetilde{\rho}_{1,2}L}\right)^j
\lim_{L\to\infty}\frac{1}{L}\sum_{j=0}^{L-1}
\left(1-\frac{2C-2C_1}{\widetilde{\rho}_{1,2}L}\right)^j\\
\end{aligned}
\end{equation}
On the other hand, the right-hand side of \eqref{SCP7} can be
rewritten
\begin{equation}\label{SCP15}
\begin{aligned}
&\psi(C) \lim_{L\to\infty}\frac{1}{L}\sum_{j=0}^{L-1}
\left(1-\frac{2C}{\widetilde{\rho}_{1,2}L}\right)^j\\
&=\psi(C) \lim_{L\to\infty}\frac{1}{L}\sum_{j=0}^{L-1}
\left(1-\frac{2C_1}{\widetilde{\rho}_{1,2}L}\right)^j
\left(1-\frac{2C-2C_1}{\widetilde{\rho}_{1,2}L}\right)^j\\
&=\psi(C)\lim_{L\to\infty}\frac{1}{L}\sum_{j=0}^{L-1}
\left(1-\frac{2C_1}{\widetilde{\rho}_{1,2}L}\right)^j
\lim_{L\to\infty}\frac{1}{L}\sum_{j=0}^{L-1}
\left(1-\frac{2C-2C_1}{\widetilde{\rho}_{1,2}L}\right)^j.
\end{aligned}
\end{equation}
The last equality in \eqref{SCP15} is an application of
\eqref{SCP14.1}. From \eqref{SCP14} and \eqref{SCP15} we finally
obtain $\psi(C_1)\geq\psi(C)$ for any positive $C_1\leq C$. The
first statement of Lemma \ref{lem6} is proved. The proof of the
second statement of this lemma is similar.
\end{proof}

We are ready now to formulate and prove a main theorem on optimal
control of the dam model considered in the present paper.

\begin{thm}\label{thm3}
The solution to the control problem is $\rho_1=1$ if and only if
the minimum of the both functionals $J^{upper}$ and $J^{lower}$ is
achieved for $C=0$. In this case the minimum of the functional $J$
is given by \eqref{SP21}. Otherwise, there can be one of the
following two cases for the solution to the control problem.

(1) If the minimum of $J^{upper}$ is achieved for $C=0$, then the
minimum of $J^{lower}$ must be achieved for some positive value
$C=\underline{C}$. Then the solution to the control problem is
achieved for $\rho_1=1-\delta$, where $\delta>0$ vanishes such
that $\delta L\to\underline{C}$ as $L\to\infty$.

(2) If the minimum of $J^{lower}$ is achieved for $C=0$, then the
minimum of $J^{upper}$ must be achieved for some positive value
$C=\overline{C}$. Then the solution to the control problem is
achieved for $\rho_1=1+\delta$, where $\delta>0$ vanishes such
that $\delta L\to\overline{C}$ as $L\to\infty$.

The minimum of the functionals $J^{upper}$ or $J^{lower}$ is
defined by their differentiating and then equating these
derivatives to zero.

\end{thm}

\begin{proof}
Note first, that there is a unique solution of the control problem
considered in this paper. Indeed, in the case where the water
costs are not taken into account, the existence of a unique
solution of the control problem follows from the main result of
\cite{Abramov 2007}. In the case of the dam model of this paper
the only difference is in presence of the functions $c^{upper}$
and $c^{lower}$ in the functionals $J^{upper}$ and $J^{lower}$
correspondingly. According to Lemma \ref{lem5} $c^{upper}=\psi(C)$
and $c^{lower}=\eta(C)$, and according to Lemma \ref{lem6} the
function $\psi(C)$ is non-increasing in $C$, while the function
$\eta(C)$ is a non-decreasing in $C$, and $\psi(0)=\eta(0)=c^*$.
Therefore, there is a unique solution of the control problem
considered in the present paper as well.

In the case where the both minima of $J^{upper}$ and $J^{lower}$
are achieved in $C=0$, then each of these minima is equal to the
right-hand side of \eqref{SP21}, and the minimum of the objective
function $J$ is achieved for $\rho_1=1$.

If the minimum of $J^{lower}$ is achieved for $C=\underline{C}$,
then $J^{lower}$ is less than the right-hand side of \eqref{SP21},
and, because of uniqueness of the solution, the minimum of
$J^{upper}$ must be achieved for $C=0$. In this case the minimum
of the objective function $J$ is achieved for $\rho_1=1-\delta$,
$\delta>0$ vanishes as $L\to\infty$, and
$L\delta\to\underline{C}$.

In the opposite case, if the minimum of $J^{upper}$ is achieved
for $C=\overline{C}$, then $J^{upper}$ is less than the right-hand
side of \eqref{SP21}, and the minimum of $J^{lower}$ must be
achieved for $C=0$. In this case the minimum of objective function
$J$ is achieved for $\rho_1=1+\delta$, $\delta>0$ vanishes as
$L\to\infty$, and $L\delta\to\overline{C}$.
\end{proof}

From Theorem \ref{thm3} we have the following evident property of
the optimal control.

\begin{cor}\label{cor2}
The solution of the control problem can be $\rho_1=1$ only in the
case $j_1\leq j_2\frac{\rho_2}{1-\rho_2}$. Specifically, the
equality is achieved only for $c_i\equiv c$, $i=1,2,\ldots,L$,
where $c$ is an any positive constant.
\end{cor}

Although Corollary \ref{cor2} provides a partial answer to the
question 2 posed in the introduction, the answer is not
satisfactory, because it is an evident extension of the result of
\cite{Abramov 2007}. A more constructive answer the question 2 of
the introduction is obtained for the special case considered in
the next section.

\section{Example of linear costs}\label{Examples}

In this section we study an example related to the case of linear
costs.

Assume that $c_1$ and $c_L<c_1$ are given.  Then the assumption
that the costs are linear means, that
\begin{equation}\label{E1}
c_i = c_1 - \frac{i-1}{L-1}(c_1-c_L), \ \ i=1,2,\ldots, L.
\end{equation}
It is assumed that as $L$ is changed, the costs are recalculated
as follows. The first and last values of the cost $c_1$ and $c_L$
remains the same. Other costs in the intermediate points are
recalculated according to \eqref{E1}.

Therefore, to avoid confusing with the appearance of the index $L$
for the fixed (unchangeable) values of cost  $c_1$ and $c_L$, we
use the other notation: $c_1=\overline{c}$ and
$c_L=\underline{c}$. Then \eqref{E1} has the form
\begin{equation}\label{E2}
c_i = \overline{c} - \frac{i-1}{L-1}(\overline{c}-\underline{c}),
\ \ i=1,2,\ldots, L.
\end{equation}

In the following we shall also use the inverse form of \eqref{E2}.
Namely,
\begin{equation}\label{E3}
c_{L-i}=\underline{c}+\frac{i}{L-1}(\overline{c}-\underline{c}), \
\ i=0,1,\ldots, L-1.
\end{equation}

Apparently,
\begin{equation}\label{E4}
c^*=\frac{\overline{c}+\underline{c}}{2}.
\end{equation}
For $c^{upper}$ we have
\begin{equation}\label{E5}
\begin{aligned}
c^{upper}=\psi(C)&=\lim_{L\to\infty}\frac{\sum_{j=0}^{L-1}
\left(\underline{c}+\dfrac{j}{L-1}(\overline{c}-\underline{c})\right)
\left(1-\dfrac{2C}{\widetilde{\rho}_{1,2}L}\right)^j}{\sum_{j=0}^{L-1}
\left(1-\dfrac{2C}{\widetilde{\rho}_{1,2}L}\right)^j}\\
&=\underline{c}+(\overline{c}-\underline{c})\lim_{L\to\infty}\frac{1}{L-1}\cdot
\frac{\sum_{j=0}^{L-1} j
\left(1-\dfrac{2C}{\widetilde{\rho}_{1,2}L}\right)^j}{\sum_{j=0}^{L-1}
\left(1-\dfrac{2C}{\widetilde{\rho}_{1,2}L}\right)^j}\\
&=\underline{c}+(\overline{c}-\underline{c})\frac{\widetilde{\rho}_{1,2}}{2C}\cdot
\frac{-\dfrac{2C}{\widetilde{\rho}_{1,2}}+
\mbox{e}^{2C/\widetilde{\rho}_{1,2}}-1}{\mbox{e}^{2C/\widetilde{\rho}_{1,2}}-1}.
\end{aligned}
\end{equation}
For example, as $C$ converges to zero in \eqref{E5}, then
$c^{upper}$ converges to
$\underline{c}+\frac{1}{2}(\overline{c}-\underline{c})=c^*$. This
is in agreement with the statement of Proposition \ref{prop1}.

In turn, for $c^{lower}$ we have
\begin{equation}\label{E6}
\begin{aligned}
c^{lower}=\eta(C)&=\lim_{L\to\infty}\frac{\sum_{j=0}^{L-1}
\left(\underline{c}+\dfrac{j}{L-1}(\overline{c}-\underline{c})\right)
\left(1+\dfrac{2C}{\widetilde{\rho}_{1,2}L}\right)^j}{\sum_{j=0}^{L-1}
\left(1+\dfrac{2C}{\widetilde{\rho}_{1,2}L}\right)^j}\\
&=\underline{c}+(\overline{c}-\underline{c})\lim_{L\to\infty}\frac{1}{L-1}\cdot
\frac{\sum_{j=0}^{L-1} j
\left(1+\dfrac{2C}{\widetilde{\rho}_{1,2}L}\right)^j}{\sum_{j=0}^{L-1}
\left(1+\dfrac{2C}{\widetilde{\rho}_{1,2}L}\right)^j}\\
&=\underline{c}+(\overline{c}-\underline{c})
\frac{\widetilde{\rho}_{1,2}}{2C}\cdot\frac{\dfrac{2C}{\widetilde{\rho}_{1,2}}-
1+\mbox{e}^{-2C/\widetilde{\rho}_{1,2}}}{1-\mbox{e}^{-2C/\widetilde{\rho}_{1,2}}}.
\end{aligned}
\end{equation}
Again, as $C$ converges to zero in \eqref{E6}, then $c^{lower}$
converges to
$\underline{c}+\frac{1}{2}(\overline{c}-\underline{c})=c^*$.
Again, we arrive to an agreement with the statement of Proposition
\ref{prop1}.

\smallskip
Let us now discuss the question 2 posed in the introduction. We
cannot give the explicit solution because the calculations are
very routine and cumbersome. However, we explain the way of the
solution of this problem and find a numerical result.

Following Corollary \ref{cor2}, take first
$j_1=j_2\frac{\rho_2}{1-\rho_2}$. Clearly, that for these relation
between parameters $j_1$ and $j_2$ the minimum of $J^{lower}$ must
be achieved for $C=0$, while the minimum of $J^{upper}$ must be
achieved for a positive $C$. Now, keeping $j_1$ fixed assume that
$j_2$ increases. Then, the problem is to find the value for
parameter $j_2$ such that the value $C$ corresponding to the
minimization problem of $J^{upper}$ reaches the point 0.

In our example we take $j_1=1$, $\rho_2=\frac{1}{2}$,
$\underline{c}=1$, $\overline{c}=2$, $\widetilde{\rho}_{1,2}=1$.
In the table below we outline some values $j_2$ and the
corresponding value $C$ for optimal solution of functional
$J^{upper}$. It is seen from the table that the optimal value is
achieved in the case $j_2\approx1.34$. Therefore, in the present
example $j_1=1$ and $j_2\approx1.34$ lead to the optimal solution
$\rho_1=1$.

\begin{table}
    \begin{center}
        \begin{tabular}{c|c}\hline
         Parameter & Argument of optimal value\\
        $j_2$ & $C$\\
        \hline
         1.06 & 0.200\\
         1.07 & 0.190\\
         1.08 & 0.182\\
         1.09 & 0.174\\
         1.10 & 0.165\\
         1.11 & 0.156\\
         1.12 & 0.149\\
         1.13 & 0.140\\
         1.14 & 0.134\\
         1.15 & 0.126\\
         1.16 & 0.120\\
        1.17 & 0.112\\
        1.18 & 0.104\\
        1.19 & 0.096\\
        1.20 & 0.090\\
        1.25 & 0.055\\
        1.30 & 0.022\\
        1.33 & 0.010\\
        1.34 & 0\\
        \hline
        \end{tabular}

        \medskip
        \caption{The values of parameter $j_2$ and corresponding
        arguments of optimal value $C$}
    \end{center}
\end{table}

\section*{Acknowledgement}
The author acknowledges with thanks the support of the Australian
Research Council, grant \# DP0771338.

% ----------------------------------------------------------------
%\bibliographystyle{amsplain}
%\bibliography{}

\end{document}